\documentclass[12pt]{article}

\usepackage{geometry}
\usepackage[affil-it]{authblk}

\usepackage{graphicx} 
\usepackage{amsmath,amsfonts}
\usepackage{xcolor}
\usepackage{booktabs}
\usepackage{hyperref}

\usepackage{lipsum}

\usepackage{authblk}

\usepackage[most]{tcolorbox}
\usepackage{lmodern}

\usepackage{fancyhdr}

 \AtEndDocument{\label{lastpage}}
\usepackage{float}
\newtcolorbox[auto counter]{mybox}[2][]{%
breakable,
enhanced,
sharp corners,
colback=white,
fonttitle=\bfseries,
title= #2,
enlarge bottom at break by=5mm,
enlarge top at break by=5mm,
overlay first={%
    \draw[black, line width=0.5mm](frame.south west)--(frame.south east);
    },
overlay middle={%
    \draw[black, line width=0.5mm](frame.south west)--(frame.south east);
    \draw[black, line width=0.5mm](frame.north west)--(frame.north east);
    \node[anchor=north east] at (frame.south east) {continued on next page};
    \node[anchor=south west] at (frame.north west) {continued from next page};
    },
overlay last={%
    \draw[black, line width=0.5mm](frame.north west)--(frame.north east);},
#1
}

\title{Revisiting the Sleeping Beauty problem}

\date{\today}

\begin{document}



\author[1]{Piva P. S.
\thanks{Electronic address: \texttt{pspiva1@sheffield.ac.uk}; Corresponding author}}
\affil{Department of Mechanical Engineering\\ The University of Sheffield\\ UK}

\author[2]{Ruffolo G.
\thanks{Electronic address: \texttt{gabrielruffolo.gr@gmail.com}}}
\affil{``Gleb Wataghin'' Institute of Physics\\ Campinas State University\\ Brazil}

\maketitle




\begin{abstract}
    The Sleeping Beauty problem is a probability riddle with no definite solution for more than two decades and its solution is of great interest in many fields of knowledge. There are two main competing solutions to the problem: the halfer approach, and the thirder approach. The main reason for disagreement in the literature is connected to the use of different probability spaces to represent the same probabilistic riddle. In this work, we analyse the problem from a mathematical perspective, identifying probability distributions induced directly from the thought experiment's rules. The precise choices of probability spaces provide both halfer and thirder solutions to the problem. To try and decide on which approach to follow, a criterion involving the information available to Sleeping Beauty is proposed. 
\end{abstract}

\section{Introduction}
\label{sec:intro}

Science deals with uncertainty and errors all the time, and one of the most useful tools to deal with this is \emph{probability theory}. Despite being an old field of mathematics and having solid foundations, there are still puzzling problems that defy common knowledge in how to use its tools precisely. Some examples are the \emph{Monty Hall problem} \cite{Monty_Hall_a,Monty_Hall_b} and the \emph{Birthday problem} \cite{birthday_1,birthday_2}. What we learn from these past problems, however, is that just intuitive reasoning is not enough to deal with probabilities, so we must be careful about the definition of the sample space, and how to correctly probabilities for its events. 

One of these puzzling riddles without clear consent for a solution is the Sleeping Beauty (SB) problem, brought to light by Elga \cite{Elga_origins,DrEvil}. In this problem, a thought experiment takes place involving a person (generically described as Sleeping Beauty or SB) who is put to sleep by some experimentalists on a Sunday. A fair coin is tossed by the experimentalists. If the outcome is Heads, SB is awakened on Monday and put back to sleep until Wednesday. If the outcome is Tails, she is awakened on Monday and put back to sleep but in this case with all her memories from Monday being erased. Then she is awakened again on Tuesday and put to sleep a final time. Finally, for any landing of the coin, she awakes on Wednesday and the experiment ends. SB cannot wake up by herself during the experiment, and she has no information about which day she woke up until Wednesday. The question of this riddle is: \emph{during an awakening on Monday or Tuesday, what is the credence of SB about the coin to land Heads?}

As stated by Elga himself \cite{Elga_origins}, the first answer one can think of is that her credence is $1/2$, since we know \emph{ab initio} that the coin is fair and SB does not seem to acquire new information about the coin during the experiment to change the probabilities. This position is defended by Lewis in his answer to Elga \cite{Lewis}, and all people who follow the same reasoning are known as \emph{Halfers}. The Halfer position has its grounds on the principle that the \emph{degree of belief} one sustains about some fact is not changed but with new information is acquired about that fact. This, in turn, is strongly related to the principles of \emph{Baysianism} \cite{Fraassen}.

Contrary to the Halfer position, Elga introduces the \emph{Thirder position} arguing that the correct answer to the problem is that the credence of SB on the result being Heads should be $1/3$ \cite{Elga_origins}. The argument is based on \emph{The Principle of Indifference}, which states that the credence about the events ``Heads \emph{and} awake on Monday'', ``Tails \emph{and} awake on Monday'' and ``Tails \emph{and} awake on Tuesday'' should be equal and, by the axioms of Kolmogorov \cite{K_axioms}, the credence of SB on the coin landing Heads is 1/3. 

Due to its implication in many different areas of knowledge 
\cite{many_areas}, this problem was extensively studied in the literature from both philosophical and mathematical perspectives \cite{Elga_origins, DrEvil, Lewis, Nick, Luna, Cozic, Groisman,Rosenthal, Barbosa}. Even hybrid models \cite{Nick} or the possibility of both positions being correct were already discussed \cite{Groisman}.

What we show in this work is that, even using the same mathematical framework developed in section \ref{sec:non-uniform}, we find (from a mathematical perspective) that both halfer and thirder positions provide possible solutions in sections \ref{sec:halfer} and \ref{sec:thirder}. However, when defining the appropriate probability spaces for the each approach, we identify the key assumption leading to each of the interpretation of SB's riddle. Then, in section \ref{sec:issue_halfer}, we introduce a new character in the experiment, who has precisely the same information as SB and, from Aumann's theorem of common knowledge \cite{Aumann}, we conclude his credence should be the same as SB's.


\section{Non-uniform probability spaces}
\label{sec:non-uniform}

A non-uniform probability space with sample space $\mathcal O$ has a probability distribution which assigns different values for different events in the sample space. In such a case, there is no combinatorial formula to calculate the probability of a specific event occurring, and a particular ad hoc value needs to be assigned to each event $p(A)$ for $A \in \mathcal O$.

However, when there is a function $f: \Omega \to \mathcal O$ that takes elements from a finite uniform sample space $\Omega$ into some of the events $A \in \mathcal O$, there is a probability distribution induced by $f$, which reads
\begin{equation}
    \label{def:induced_prob}
    p(A) := \frac{\# f^{-1}(A)}{\# \Omega},
\end{equation}
where $\#$ denoted the cardinality\footnote{Equivalent to the number of elements of the set, for the case of finite sets.}, and $f^{-1}$ is the pre-image of $f$. If $\# \mathcal O > \# f(\Omega)$ , then $p(A) := 0$ for any $A \notin f(\Omega)$.

The probability distribution is not defined only for elements in the sample space but over a bigger set called the event space $\mathfrak F (\mathcal O)$. Then, we define $p: \mathfrak F (\mathcal O) \to [0,1]$ with the usual extension
\begin{equation}
    \label{def:prob_event_space}
    p(A) := \frac{\# \{ X \in \Omega | X \subseteq f^{-1} (A) \}}{\# \Omega},
\end{equation}
for all $A \in \mathrm P(f(\Omega))$, where $\mathrm P(\, . \,)$ denotes the power set. If $A \notin \mathrm P(f(\Omega))$ we define
\begin{equation}
    \label{def:prob_out_induced}
    p(A) := p \left( A \setminus \overline{f(\Omega)} \right),
\end{equation}
where $\overline{f(\Omega)}$ denotes the complement of the image of $f$.

If $f$ is also a defining property of the sample space $\mathcal O$, then \eqref{def:induced_prob} is the only possible probability distribution for the sample space by definition, in the sense that any other ad hoc probability distribution chosen will not share such property.

\begin{mybox}{Example 1}
    
One example from introductory probability books \cite{Ross} is the sum of two 6-sided fair dice, which we represent by $\mathcal S$.

The finite uniform sample space that represents the action of rolling two 6-sided fair dice, which reads
\begin{equation}
    \label{def:sample_2_dice}
    \Omega_D := \left\{ (i, j) \in \mathbb N^2 | 1\leq i, j \leq 6 \right\}
\end{equation}
The function which induces a probability on $\mathcal S$ is $s: \Omega_D \to \mathcal S$, given by 
\begin{equation}
    \label{def:sum_function}
    s(i,j) := i + j.
\end{equation}
Because $s(\Omega_D)$ covers only the integers from 2 up to 12, we have that $p(s = N) = 0$ for any $N < 2$ or $N > 12$. Four explicit values of probability are calculated below using \eqref{def:induced_prob}, \eqref{def:prob_event_space} and \eqref{def:prob_out_induced}.
\begin{equation}
    \notag
    \begin{aligned}
        &p(s = 2) = \frac{\# s^{-1}(2)}{\# \Omega_D} = \frac{\#\{ (1,1) \}}{36} = \frac{1}{36},
\\
        &p(s = 6) = \frac{\# s^{-1}(6)}{\# \Omega_D} = \frac{\#\{ (1,5), (2,4), (3,3), (4,2), (5,1) \}}{36} = \frac{5}{36},
\\
        &p(s = 2 \text{ or } 6) = \frac{\# \{ X \in \Omega| X \subseteq s^{-1}(\{2,6\}) \}}{\# \Omega_D} = \frac{6}{36} = \frac{1}{6},
\\
        &p(s < 3) = p ((s < 3) \setminus (s < 2)) = p (s = 2) = \frac{1}{36}.
    \end{aligned}
\end{equation}

Any other choice of probability distribution associated with $\mathcal S$ would represent a different probabilistic space than the one for the sum of two fair dice. Either the dice would not be fair, or the sum rule $s$ would not be respected.
    
\end{mybox}

\subsection{Joint and conditional probability}
\label{subsec:conditional}

Given the definition \eqref{def:induced_prob}, a natural question is whether it is possible or not to define conditional probabilities concerning events in both $\Omega$ and $\mathcal O$. To achieve this goal, we first define the joint probability of two events, one in each sample space.

{\em Definition:} let $\Omega$ be a finite uniform sample space and $\mathcal O$ any sample space. Let $f: \Omega \to \mathcal O$ be a function which induces the probability distribution assigned to $\mathcal O$. We define the joint probability of an event $A \in \mathfrak F(\mathcal O)$ and $B \in \mathrm P(\Omega)$ as
\begin{equation}
    \label{def:joint_prob}
    p(A \, \& \, B) := \frac{\#\{X \in \Omega | X \subseteq \left( f^{-1}(A) \cap B \right) \}}{\#\Omega},
\end{equation}
which gives us $p(A \, \& \, B) \leq p(A)$ or $p(B)$.

For the case where $\mathcal O$ is enumerable, we have that the sum of $p(A \, \& \, B)$ for all $A \in \Omega$ and $B \in \mathcal O$ is 1, by definition. This means we can define the product sample space $\mathcal P = \Omega \times \mathcal O$, whose probability distribution is given by \eqref{def:joint_prob}.

\begin{mybox}{Example 2}
    
With the same setup for Example 1, we define the sample space $\Omega_D \times \mathcal S$, and assign the joint probability induced by \eqref{def:sum_function} as defined in \eqref{def:joint_prob}. For an explicit calculation, we choose the event $(s = 6) \, \& \, B$, where $B \in \mathrm P(\Omega_D)$ s.t. $B =$ ``At least one of the dice scored 1''. The result reads
\begin{equation}
    \notag
    p ((s = 6) \, \& \, B) = \frac{\#\{ X \in \Omega|X \subseteq \left( s^{-1}(6) \cap B \right) \}}{\# \Omega_D} = \frac{\#\{ (1,5), (5,1) \}}{36} = \frac{1}{18}.
\end{equation}

\end{mybox}

Using the probability distributions associated with all sample spaces $\Omega$, $\mathcal O$ and $\mathcal P$, we can define the conditional probabilities as
\begin{equation}
    \label{def:cond_prob}
    \begin{aligned}
        p(A|B) := \frac{p(A \, \& \, B)}{p(B)},
\\
        p(B|A) := \frac{p(A \, \& \, B)}{p(A)},
    \end{aligned}
\end{equation}
where we have used the uniform probability distribution of $\Omega$ to compute $p(B)$.

\begin{mybox}{Example 3}
    
We recall Examples 1 and 2, and we consider the events $s = 6$ and $B =$ ``At least one of the dice scored 1''. We use the definition \eqref{def:cond_prob} to explicitly compute $p(s = 6|B)$ below.
\begin{equation}
    \notag
    p (s = 6|B) = \frac{p ((s = 6) \, \& \, B)}{p(B)} = \frac{1}{18} \times \frac{36}{11} = \frac{2}{11}.
\end{equation}

\end{mybox}

\subsection{Credence functions}
\label{subsec:credence}

In this section, we discuss how to compute credence with a sample space representing the roll of two fair dice and its corresponding sum, as defined in Example 1 by \eqref{def:sample_2_dice} and \eqref{def:sum_function}.

We gather two fictitious characters, named Alice and Bob, and tell them we will be rolling two fair 6-sided dice. We roll the dice and take the record of both scores $d_1 = 5$ and $d_2 = j$ with $1 \leq j \leq 6$, but we only tell the first score $d_1 = 5$ to Alice and give Bob no information about either scores or the knowledge of Alice about $d_1$. Then, we ask them which is their credence on the result of the total score $s = d_1 + d_2$ for all possible values of $s$.

Both Alice and Bob know probability theory and all the definitions presented so far in the text. Alice's credence is given by the probability of $s$ conditional to the event $d_1 = 5$, which reads $C_{A}(s) = p(s|d_1 = 5)$. Bob's credence is given by simply $C_{B}(s) = p(s)$, because he does not have any additional information about the scores. Then, we gather both our helpers in the same room and ask Alice if she believes the total score may be $s = 4$. After Alice answers that $s = 4$ is impossible, Bob concludes she has more information than him about either $d_1$ or $d_2$ and his credence gets updated to $C'_{B} = p(s|s \neq 4)$. Finally, we tell Bob we had informed Alice about the result of $d_1$, so he infers that $d_1 \geq 4$, and his credence gets updated once again, this time to $C''_{B}(s) = p(s|d_1 \geq 4)$. All credence functions of both Alice and Bob at all times of the experiment are represented in Table \ref{tab:credences_1} below.
\begin{table}[h]
\centering
\begin{tabular}{@{}cccccccccccc@{}}
$s$ & 2 & 3 & 4 & 5 & 6 & 7 & 8 & 9 & 10 & 11 & 12 \\ \midrule
$C_A(s)$ & 0 & 0 & 0 & 0 & $\frac{1}{6}$ & $\frac{1}{6}$ & $\frac{1}{6}$ & $\frac{1}{6}$ & $\frac{1}{6}$ & $\frac{1}{6}$ & 0 \\ \midrule
$C_B(s)$  & $\frac{1}{36}$ & $\frac{1}{18}$ & $\frac{1}{12}$ & $\frac{1}{9}$ & $\frac{5}{36}$ & $\frac{1}{6}$ & $\frac{5}{36}$ & $\frac{1}{9}$ & $\frac{1}{12}$ & $\frac{1}{18}$ & $\frac{1}{36}$ \\ \midrule
$C'_B(s)$ & $\frac{1}{33}$ & $\frac{2}{33}$ & 0 & $\frac{4}{33}$ & $\frac{5}{33}$ & $\frac{6}{33}$ & $\frac{5}{33}$ & $\frac{4}{33}$ & $\frac{1}{11}$ & $\frac{2}{33}$ & $\frac{1}{33}$ \\ \midrule
$C''_B(s)$ & 0 & 0 & 0 & $\frac{1}{18}$ & $\frac{1}{9}$ & $\frac{1}{6}$ & $\frac{1}{6}$ & $\frac{1}{6}$ & $\frac{1}{6}$ & $\frac{1}{9}$ & $\frac{1}{18}$
\end{tabular}
\caption{Table of credence functions of Alice and Bob for the sum of two dice. Alice's credence is unaltered throughout the experiment, represented by $C_A$. Bob's credence gets updated twice and is represented by $C_B$, $C_B'$ and $C_B''$.}
\label{tab:credences_1}
\end{table}

A common mistake when computing credence is to forget about the complex structure of the sample space $\mathcal S$, and redefine a new probability distribution for Bob each time he learns something about $d_1$ and $d_2$. For example, one could assume that all the proportions between probabilities of $\mathcal S$ are kept the same, apart from the events which Bob learns to be impossible to happen, then the credence could be computed like in Table \ref{tab:credences_2} below. This ad hoc procedure leads to a different result for the updated credence of Bob in Table \ref{tab:credences_1} because it does not take into consideration the correlations between different values of $s$ induced by the function $s$ itself.
\begin{table}[h]
\centering
\begin{tabular}{@{}cccccccccccc@{}}
$s$ & 2 & 3 & 4 & 5 & 6 & 7 & 8 & 9 & 10 & 11 & 12 \\ \midrule
$C_B(s)$  & $\frac{1}{36}$ & $\frac{1}{18}$ & $\frac{1}{12}$ & $\frac{1}{9}$ & $\frac{5}{36}$ & $\frac{1}{6}$ & $\frac{5}{36}$ & $\frac{1}{9}$ & $\frac{1}{12}$ & $\frac{1}{18}$ & $\frac{1}{36}$ \\ \midrule
$C'_B(s)$ & $\frac{1}{33}$ & $\frac{2}{33}$ & 0 & $\frac{4}{33}$ & $\frac{5}{33}$ & $\frac{6}{33}$ & $\frac{5}{33}$ & $\frac{4}{33}$ & $\frac{1}{11}$ & $\frac{2}{33}$ & $\frac{1}{33}$ \\ \midrule
$C''_B(s)$ & 0 & 0 & 0 & $\frac{2}{15}$ & $\frac{1}{6}$ & $\frac{1}{5}$ & $\frac{1}{6}$ & $\frac{2}{15}$ & $\frac{1}{10}$ & $\frac{1}{15}$ & $\frac{1}{30}$
\end{tabular}
\caption{Table of credence functions for Bob using an ad hoc sample space at each update. Once again, Bob's credence is updated twice and is represented by $C_B$, $C_B'$ and $C_B''$.}
\label{tab:credences_2}
\end{table}

\section{Solution from halfer's perspective}
\label{sec:halfer}

To understand the Sleeping Beauty (SB) problem from the halfer's perspective, we define the uniform sample space representing the toss of a fair coin as
\begin{equation}
    \notag
    \Omega_C := \{ H, T \},
\end{equation}
where $H$ and $T$ are the two possible landings of ``Heads'' or ``Tails'' respectively. We also define the following sample space
\begin{equation}
    \label{eq:sample_example}
    \mathcal B = \{ (0,0), (1,0), (0,1), (1,1)\},
\end{equation}
and attribute the following meaning for each of its elements.
\begin{equation}
    \label{eq:days_association}
    \begin{aligned}
        & (0, 0) \to \text{``SB is not awakened on Monday neither on Tuesday''};
\\
        & (1, 0) \to \text{``SB is awakened on Monday, but not on Tuesday''};
\\
        & (0, 1) \to \text{``SB is awakened on Tuesday, but not on Monday''};
\\
        & (1, 1) \to \text{``SB is awakened on Monday and Tuesday''},
    \end{aligned}
\end{equation}
At last, we define the function $g:\Omega_C \to \mathcal B \times \mathcal B$
\begin{equation}
    \label{def:function_coin}
    \begin{aligned}
        g(H) := (1,0),
\\
        g(T) := (1,1).
    \end{aligned}
\end{equation}
which describes the rules of the experiment.

Using \eqref{def:joint_prob} we compute the joint probabilities induced in the product space $\Omega_C \times \mathcal B$, given by
\begin{equation}
    \label{eq:probs_halfer}
    \begin{aligned}
        & p(H \, \& \, (0,0)) = 0;
\\
        & p(H \, \& \, (1,0)) = 1/2;
\\
        & p(T \, \& \, (0,1)) = 0;
\\
        & p(T \, \& \, (1,1)) = 1/2.
    \end{aligned}
\end{equation}

In possession of the joint probabilities \eqref{eq:probs_halfer} above, we continue the discussion below, computing in detail both credence functions of the experimentalists and SB herself.

\subsection{Credence of experimentalists}
\label{subsec:halfer_experimentalists}

We start by defining the projections of the binary space of waking up SB or not on Monday, $h_M:\mathcal B \times \mathcal B \to \mathcal B_M$, and Tuesday, $h_T:\mathcal B \times \mathcal B \to \mathcal B_T$, as
\begin{equation}
    \notag
    \left\{
    \begin{aligned}
        & h_M(0, 0) := 0;
\\    
        & h_M(1, 0) := 1;
\\
        & h_M(0, 1) := 0;
\\
        & h_M(1, 1) := 1;
    \end{aligned}
    \right.
\end{equation}
\begin{equation}
    \notag
    \left\{
    \begin{aligned}
        & h_T(0, 0) := 0;
\\    
        & h_T(1, 0) := 0;
\\
        & h_T(0, 1) := 1;
\\
        & h_T(1, 1) := 1,
    \end{aligned}
    \right.
\end{equation}
where $\mathcal B_M$ and $\mathcal B_T$ are binary sample spaces, and we attribute the following meaning to each of their events
\begin{equation}
    \label{eq:days_association_2}
    \begin{aligned}
        & A_M = 1 \to \text{``SB is awakened on Monday''};
\\
        & A_M = 0 \to \text{``SB is not awakened on Monday''};
\\
        & A_T = 1 \to \text{``SB is awakened on Tuesday''};
\\
        & A_T = 0 \to \text{``SB is not awakened on Tuesday''},
    \end{aligned}
\end{equation}
for $A_M \in \mathcal B_M$ and $A_T \in \mathcal B_T$. This allows us to compute the probabilities induced by the composition $g \circ h_M: \Omega_c \to \mathcal B_M$ and $g \circ h_T: \Omega_c \to \mathcal B_T$, which leads to
\begin{equation}
    \notag
    \label{eq:M_probs}
    \left\{
    \begin{aligned}
        & p_M(A_M = 1) = 1;
\\
        & p_M(A_M = 0) = 0.
    \end{aligned}
    \right.
\end{equation}
\begin{equation}
    \notag
    \label{eq:T_probs}
    \left\{
    \begin{aligned}
        & p_T(A_T = 1) = \frac{1}{2};
\\
        & p_T(A_T = 0) = \frac{1}{2}.
    \end{aligned}
    \right.
\end{equation}

Then, using definition \eqref{def:joint_prob} to compute the respective joint probabilities induced by $g \circ h_M$ and $g \circ h_T$, we have
\begin{equation}
    \label{eq:M_joint_probs}
    \left\{
    \begin{aligned}
        & p_M(H \, \& \, (A_M = 1)) = 1/2;
\\
        & p_M(H \, \& \, (A_M = 0)) = 0;
\\
        & p_M(T \, \& \, (A_M = 1)) = 1/2;
\\
        & p_M(T \, \& \, (A_M = 0)) = 0;
    \end{aligned}
    \right.
\end{equation}
\begin{equation}
    \label{eq:T_joint_probs}
    \left\{
    \begin{aligned}
        & p_T(H \, \& \, (A_T = 1)) = 0;
\\
        & p_T(H \, \& \, (A_T = 0)) = 1/2;
\\
        & p_T(T \, \& \, (A_T = 1)) = 1/2;
\\
        & p_T(T \, \& \, (A_T = 0)) = 0.
    \end{aligned}
    \right.
\end{equation}
Both probability densities \eqref{eq:M_joint_probs} and \eqref{eq:T_joint_probs} are assigned to both product spaces  $\mathcal C_M := \Omega_C \times \mathcal B_M$ and $\mathcal C_T := \Omega_C \times \mathcal B_T$ respectively.

Each of the probabilities above will represent the credence functions of all possible outcomes for the experimentalists before tossing the coin. Even if they toss the coin on Tuesday, right before knowing if they should or not wake up SB, these credences would remain the same until the tossing. It is important to notice that \eqref{eq:M_joint_probs} and \eqref{eq:T_joint_probs} \emph{agree} with the Principle of Indifference because we have that
\begin{equation}
    p_M(H \, \& \, (A_M = 1)) = p_M(T \, \& \, (A_M = 1)) = p_T(T \, \& \, (A_T = 1)).
\end{equation}

\subsection{Credence of Sleeping Beauty}
\label{subsec:halfer_SB}

Because SB knows her memory will be erased if her fate is to be awakened twice, she cannot know which day she woke up. She can only ask herself if she is awake or not at the moment.

Mathematically, the previous facts translate into her credence being given by the probability function induced by \eqref{def:function_coin}, but conditioned to the only event which she has access to: $W = \{(1, 0), (0, 1), (1, 1) \} \to$ ``SB is awakened on Monday or Tuesday''. Using definitions \eqref{def:prob_event_space} and \eqref{def:prob_out_induced} to compute $p(W) = 1$, which we combine with \eqref{def:joint_prob} and \eqref{def:cond_prob} to determine the famous values of credence for SB from the halfers approach
\begin{equation}
    \label{eq:credence_halfer}
    \begin{aligned}
        & p(H|W) = p(H) = 1/2,
\\
        & p(T|W) = p(T) = 1/2.
    \end{aligned}
\end{equation}

Equations \eqref{eq:credence_halfer} above tell us that SB knows she will wake at some point ($p(W) = 1$), and she knows that fact from the rules of the experiment. Then, she does not acquire any new information during an awakening and does not change her mind about the chances of the fair coin to land Heads or Tails.

\section{Solution from thirder's perspective}
\label{sec:thirder}

To understand the Sleeping Beauty (SB) problem from the thirder's perspective we define the following uniform sample space
\begin{equation}
    \label{eq:sample_space_coin_day}
    \Omega_{C,D} := \{ (H, \text{Mo}), (T, \text{Mo}), (H, \text{Tu}), (T, \text{Tu}) \},
\end{equation}
where each tuple represents the following events
\begin{equation}
    \notag
    \begin{aligned}
        & (H, \text{Mo}) \to \text{``The fair coin landed Heads and today is Monday''};
\\
        & (T, \text{Mo}) \to \text{``The fair coin landed Tails and today is Monday''};
\\
        & (H, \text{Tu}) \to \text{``The fair coin landed Heads and today is Tuesday''};
\\
        & (T, \text{Tu}) \to \text{``The fair coin landed Tails and today is Tuesday''}.
    \end{aligned}
\end{equation}

We also define the binary sample space of $\mathcal A := \{ A, S \}$ with its events meaning $A \to$ ``SB will be awakened today'' and $S \to$ ``SB will keep sleeping today''.

Then, the rules of the experiment are given by the awakening function $a: \Omega_{C,D} \to \mathcal A$, defined as
\begin{equation}
    \label{def:function_coin_day}
    \begin{aligned}
        &a(H, \text{Mo}) := A;
\\
        &a(T, \text{Mo}) := S;
\\
        &a(H, \text{Tu}) := A;
\\
        &a(T, \text{Tu}) := A.
    \end{aligned}
\end{equation}

Finally, we compute all joint probabilities for the product space $\Omega_{C, D} \times \mathcal A$ using definition \eqref{def:joint_prob}, which leads to
\begin{equation}
    \label{eq:thirder_probs}
    \begin{aligned}
        &p(H, \text{Mo} \, \& \,  A) = 1/4; \quad p(H, \text{Mo} \, \& \,  S) = 0;
\\
        &p(T, \text{Mo} \, \& \,  A) = 1/4; \quad p(T, \text{Mo} \, \& \,  S) = 0;
\\
        &p(H, \text{Tu} \, \& \,  A) = 0; \quad \ \ \  p(H, \text{Tu} \, \& \,  S) = 1/4;
\\
        &p(T, \text{Tu} \, \& \,  A) = 1/4; \quad p(T, \text{Tu} \, \& \,  S) = 0.
    \end{aligned}
\end{equation}

As done in section \eqref{sec:halfer}, we continue by using the joint probabilities \eqref{eq:thirder_probs} above to compute both credence functions of the experimentalists and SB herself.

\subsection{Credence of experimentalists}
\label{subsec:thirder_experimentalists}

Because the experimentalists can keep track of the days of the week, their credence on each outcome will be given by the following
\begin{equation}
    \label{eq:joint_probs_thirder}
    \begin{aligned}
        & p(H \, \& \, A | \text{Mo}) = 1/2; \quad p(H \, \& \, S | \text{Mo}) = 0;
\\
        & p(H \, \& \, A | \text{Mo}) = 1/2; \quad p(H \, \& \, S | \text{Mo}) = 0;
\\
        & p(T \, \& \, A | \text{Tu}) = 0; \quad \quad \  p(T \, \& \, S | \text{Tu}) = 1/2;
\\
        & p(T \, \& \, A | \text{Tu}) = 1/2; \quad \  p(T \, \& \, A | \text{Tu}) = 0,
    \end{aligned}
\end{equation}
which were calculated using definitions \eqref{def:joint_prob} and \eqref{def:cond_prob}. Notice that equations \eqref{eq:joint_probs_thirder} above agree with \eqref{eq:M_joint_probs} and \eqref{eq:T_joint_probs} from the halfers solution. This means that both approaches should be valid to compute the credences of the experimentalists.

\subsection{Credence of Sleeping Beauty}
\label{subsec:thirer_SB}

Once again, SB cannot determine which day she was awakened, and the only question she can ask herself is: \emph{am I awake right now?}

In this case, her credence is given by the probability function induced by \eqref{def:function_coin_day} conditioned to the only event $A \to$ ``SB will be awakened today''. Using definition \eqref{def:prob_event_space} we determine $p(A) = 3/4$, which combined with \eqref{def:joint_prob} and \eqref{def:cond_prob} leads to
\begin{equation}
    \label{eq:credence_thirder}
    \begin{aligned}
        & p(H|A) = 1/3 \neq p(H),
\\
        & p(T|A) = 2/3 \neq p(T).
    \end{aligned}
\end{equation}
which are the also famous values of credence for SB from the thirder's approach.

To interpret these results we need to consider that SB gains new information about the coin when she is awakened, so she updates her credence about a fair coin landing Heads from $p(H)$ to $p(H|A)$, and similarly for a Tails landing. This interpretation of \eqref{eq:credence_thirder} \emph{is consistent} with the principles of Bayesianism because $p(H, \text{Tu}|A) = 0$, and SB discards one of the four events in \eqref{eq:sample_space_coin_day} when she is awakened.

\subsection{The original thirder approach}
\label{subsec:issue_thirder}

In the original thirder approach \cite{Elga_origins}, the problem was solved by computing the following
\begin{equation}
    \label{eq:Elga_probs}
    p(H \, \& \, \text{Mo} | A) = p(T \, \& \, \text{Mo} | A) = p(T \, \& \, \text{Tu} | A) = 1/3,
\end{equation}
which agree with the probability distributions \eqref{eq:joint_probs_thirder} calculated at the beginning of this section. However, the calculations do not clearly state that the awakening is considered a random variable in the analysis, leading to the conclusion that $p(H|A) = 1/3$ directly from \eqref{eq:Elga_probs}. This approach may lead to the interpretation that $p(H) = 1/3$, which contradicts the experiment's rules because the coin would not be fair from the beginning. With the mathematical framework presented in section \ref{sec:non-uniform}, we let it clear in \eqref{eq:credence_thirder} that, even though SB knows the tossed coin is fair, her credence on a Heads landing ends up being 1/3 during an awakening. 

\section{Issue with the halfer position}
\label{sec:issue_halfer}

After computing the credence of Sleeping Beauty (SB) on both Heads and Tails landings for a fair coin, we find two contradicting results, \eqref{eq:credence_halfer} and \eqref{eq:credence_thirder}, calculated using the very same strategy introduced in section \ref{sec:non-uniform}. The main difference in both approaches is to consider the event of SB's awakening and the day it happens as separate random variables (thirder approach) or not (halfer approach). In this section, a new character is introduced into the thought experiment, so we can apply an objective criterion to decide on which is the correct approach to solve this riddle.

We invite Grumpy into the thought experiment, who plays an important role in the discussion. We keep him in a separate room with only one phone. Each day SB is awakened, we give her a mobile phone to call Grumpy and tell him she is awake. His memory of the whole day will be erased at the end of Monday. Before the experiment starts, we explain all the rules regarding SB and himself to Grumpy. Then, we ask ourselves: \emph{what is the credence of Grumpy about the coin to land Heads after receiving a call from SB on any day?}

On each day, Grumpy does not know if SB will be awakened, nor which day it is, so he must consider her \emph{awakening} as a \emph{separate random variable} than the \emph{day of the awakening}. This variable is correlated with two other random variables in the problem: the fair coin landing and which day is today.

\begin{figure}[H]
    \centering
    \includegraphics[width=0.5\textwidth]{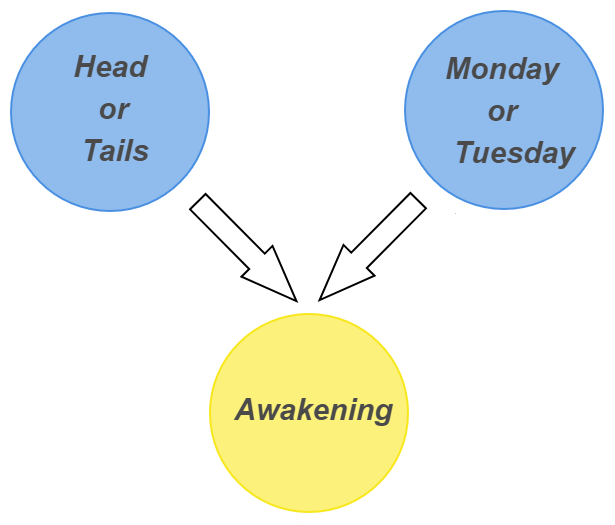}
    \caption{Causal structure of the events in the experiment.}
    \label{fig:causaldiagram}
\end{figure}

We consider the correlations implied by the rules of the experiment. The coin and the day are enough to \emph{completely determine} if the awakening of SB would happen, so in fact they have a \emph{causal structure} (see Figure \ref{fig:causaldiagram}). Then, we can say that the following correlations hold (in accordance with \eqref{def:function_coin_day})
\begin{equation}
    \label{eq:correlations}
    \begin{aligned}
        &p(A|\text{Mo},H) = 1;
\\
        &p(A|\text{Mo},T) = 1;
\\
        &p(A|\text{Tu},H) = 0;
\\
        &p(A|\text{Tu},T) = 1,
    \end{aligned}
\end{equation}

From Figure \ref{fig:causaldiagram}, we know that the result of the coin toss on Sunday is independent of the day on which Grumpy is waiting for the phone call, leading to
\begin{equation}
    \label{eq:independence_day}
    p(Day, Coin) = p(Day)p(Coin),
\end{equation}
where $Day \in \{\text{Mo},\text{Tu}\}$ and $Coin \in \{H,T\}$. We also know that the coin is fair, so $p(Coin) = 1/2$. The last thing to consider is the probability of the days alone, which Grumpy should assign $p(Day) = 1/2$ because \emph{nothing is causing it}, and we can apply the Principle of Indifference.

Using Bayes rule in each equation from \eqref{eq:correlations} and substituting \eqref{eq:independence_day}, we have
\begin{equation}
    \notag
    \begin{aligned}
        &p(A,\text{Mo},H) = 1/4;
\\
        &p(A,\text{Mo},T) = 1/4;
\\
        &p(A,\text{Tu},H) = 0;
\\
        &p(A,\text{Tu},T) = 1/4,
    \end{aligned}
\end{equation}
which agrees with the probability distributions for the thirder approach \eqref{eq:thirder_probs}. The probabilities above are enough to calculate the following marginals:
\begin{equation}
    \label{eq:marginals}
    \left\{
    \begin{aligned}
        &p(A,H) = p(A,\text{Mo},H) + p(A,\text{Tu},H) = 1/4
\\
        &p(A,T) =  p(A,\text{Mo},T) + p(A,\text{Tu},T) = 1/2
    \end{aligned}
    \right.
    \Rightarrow
    p(A) = p(A,H) + p(A,T) = 3/4.
\end{equation}

Now we can compute Grumpy's credence about the coin landing Heads, after SB called him and he knows that she is awaken. As done in sections \ref{subsec:halfer_SB} and \ref{subsec:thirer_SB}, we just need to calculate the conditional probability $p(H|A)$. Using the marginals from \eqref{eq:marginals} we conclude
\begin{equation}
    \label{eq:credence_Grumpy}
    \begin{aligned}
        p(H|A) = \frac{p(A,H)}{p(A)} = \frac{1/4}{3/4} = 1/3;
\\
        p(T|A) = \frac{p(T,H)}{p(A)} = \frac{1/2}{3/4} = 2/3,
    \end{aligned}
\end{equation}
which is in agreement with the thirder position once again.

Finally, we notice that Grumpy has \emph{the exact same information} about the experiment as SB after he answers her phone call. This allows us to use Aumann's theorem of common knowledge \cite{Aumann} to conclude that SB's credence on the landing of the fair coin is also given by \eqref{eq:credence_Grumpy}. The result is that SB also needs to consider the event of her awakening as a separate random variable from the day of the awakening.

\subsection{Issue with the ambiguity argument}
\label{subsec:ambiguity}

An attempt to reformulate the SB problem in simpler terms was given by Groisman \cite{Groisman}. He proposed we toss a fair coin and, if it lands Heads, we deposit a green ball in a box and, if lands Tails, we deposit two red balls in the same box. After many tossings, we can ask ourselves: \emph{what is our credence about the coin landing Heads?} and \emph{what is our credence about randomly picking a green ball from the box?}.

To answer these questions we need to build two sample spaces, one is the uniform space with events for coin landing Heads or Tails, given by $\Omega_C = \{ H, T \}$. The other is the space of number of balls after $N$ tossings $\mathcal A_N = \{ G, R \}$, where $0 \leq G \leq N$ and $0 \leq R \leq 2N$. The rules of this experiment are induced by $S: \left( \Omega_C\right)^N \to \mathcal A_N$, where a set to the power $N$ stands for the result of $N$ cartesian products with itself. For each tossing, $S$ adds one to $G$ if the coin lands $H$, or it adds 2 to $R$ if the coin lands tails. Before the first tossing $R = S = 0$.

With the probability spaces properly defined, we conclude that the answer to the first question is 1/2, because the coin is fair, while the answer to the second question approaches the value of 1/3 as $N \to \infty$. In \cite{Groisman}, it is argued that this example implies that both halfer and thirder positions should be valid, and the question is ambiguous if asked in the same terms as stated in section \ref{sec:intro}. However, this cannot be inferred from this example of two different questions formulated about two events belonging to \emph{completely different probability spaces} compared to the one describing the question of the SB problem.

To formulate a question about an event from a probability space with the same structure as in the SB problem, we restrict ourselves to the case of only one tossing of the coin ($N=1$) and ask, instead: \emph{what is our credence about the coin landing Heads if we pick a ball from the box regardless of its colour?} The answer to this question is undoubtedly 1/2. This example shows, once again, that the halfer position is a self-consistent solution to the SB problem from a mathematical perspective. However, this setup does not capture the nuances introduced by erasing SB's memory on Monday, so the question is still essentially different from the one regarding SB's credence.

\section{Conclusion}
\label{sec:conclusion}

In this work, we showed how to deal with the Sleeping Beauty (SB) problem in its original formulation. All the definitions and results needed for the solution of the problem are provided throughout section \ref{sec:non-uniform}. Most of the definitions depend on a function lying on a general sample space while acting on elements of a uniform sample space. This function plays an important role in the interpretation of the probability distribution assigned to the product probability space, as shown by the examples presented. The notion of credence is also discussed for a specific example of rolling two dice.

With the aid of the framework presented, we identify the appropriate probability spaces to solve the SB problem in sections \ref{sec:halfer} and \ref{sec:thirder}, which are compatible with both halfer and thirder positions. The credences of all individuals involved in the thought experiment are explicitly computed in this section. Throughout these sections, it is interesting to notice that none of the approaches harm the Principle of Indifference or the principles of Bayesianism, as usually discussed in the literature. That means another criterion is needed to decide which probability SB should choose to compute her credence.

Then, in section \ref{sec:issue_halfer}, we explain the issue of the solution provided by the halfer position by using Aumann's theorem of common knowledge \cite{Aumann}. Even though the halfer solution respects the correlations introduced by the rules of the experiments, it does not correspond to the case in which SB's memory gets erased, and she does not have information about which day she is awakened on the experiment. Under these circumstances, she needs to consider the event of her awakening as a separate random variable from the day of awakening.  However, the halfer approach is consistent from the perspective of the experimentalists and predicts the correct credences for them.

In section \ref{subsec:ambiguity}, we explain why the ambiguity argument \cite{Groisman} is also not a viable solution to the riddle. The questions formulated within this approach do not consider events from the same sample spaces with a similar structure as in SB's problem.

We call the attention to the reader that the concept of \emph{self-locating belief} \cite{Elga_origins,DrEvil,many_areas,quantum} is not needed to reach the credences from the thirder's approach. The careful analysis of which probability spaces to use and the information available to SB is enough to solve the SB problem. This suggests that the framework presented may contribute to new methodologies for analysing problems in the areas of the philosophy of science and decision theory \cite{many_areas}.


\section*{Acknowledgements}
This work has been supported by the Engineering and Physical Sciences Research Council EPSRC (Industrial CASE studentship with Johnson Matthey) and the São Paulo Research
Foundation FAPESP (Grants No. 2021/
01502-6, and No. 2021/10548-0).

\end{document}